\magnification\magstep1
 \input epsf.sty
\input  amssym.tex
\input color
%
 
 \centerline{\bf  Dynamics of Polynomial Diffeomorphisms in ${\Bbb C}^2$:}
 
 \centerline{\bf  Foliations and laminations}
 
 \bigskip
 \centerline{Eric Bedford}
 
 \bigskip
 
 Around the beginning of the 20th century, Poincar\'e considered the restricted three-body problem.  This led him to the study of surface diffeomorphisms.  These maps arose from the ODE defined by Newton's laws: either as the time-one map of the evolution, or as the first return map to a Poincar\'e section.  Poincar\'e recast the study of evolution of Newton's laws as a search for the qualitative behavior of long term solutions.  He formulated the question of dynamics as to determine the behavior of the iterates of mappings $f^n:=f\circ\cdots\circ f$ as $n\to\infty$.  Poincar\'e focused on recurrent behavior:  a point $x_0$ is recurrent if there is a sequence $n_j\to\infty$ such that $f^{n_j}(x_0)\to x_0$.  Among the recurrent points, the ones of particular interest are the periodic ones, which satisfy $f^n(x_0)=x_0$ for some $n\ne0$.  
 
 A periodic point $x_0$, with $f^n(x_0)=x_0$, is a saddle point if the tangent space at $x_0$ splits into a subspace on which the differential  $D_{x_0}f^n$ is strictly contracting and a complementary subspace on which it is strictly expanding.   In this case there is a stable manifold $W^s(x_0)$, which consists of points which approach $x_0$ in forward time.  Similarly, there is the unstable manifold $W^u(x_0)$ of points which approach $x_0$ in backward time.  The intersection points $W^s(x_0)\cap W^u(x_0)$ are said to be homoclinic.  Poincar\'e observed that when there is a homoclinic point where $W^s$ and $W^u$ have a transverse intersection, then the configuration of $W^s\cup W^u$ forms a very complicated sort of ``trellis'' or ``homoclinic tangle''.  Poincar\'e also noticed that chaotic behavior is caused by these homoclinic points, and that there is a connection between the complicated geometry and the chaotic dynamics.
  
 
A hyperbolic set is an invariant compact set on which there is uniform expansion and contraction as in the case of a saddle point.  The theory of hyperbolic maps was developed during the 1960s by Anosov, Sinai, Smale, and others.  Smale formulated his ``horseshoe map'' as a simple model of chaotic behavior: this model map is hyperbolic, and its dynamics corresponds to the random behavior of Bernoulli trials.  
The success of the theory in describing hyperbolic maps raised the question of how typical hyperbolic maps might be.  It was shown that hyperbolicity was not a generic phenomenon, and Newhouse gave examples of this in the plane.
 
For a given diffeomorphism, there is the dichotomy between the regions of stable and chaotic behavior, and the interplay between these two regions is intriguing.   A model for stable behavior would be an attracting fixed (or periodic) point $x_0$.  The basin of all points that converge to $x_0$ is an open set, and the dynamical behavior on the basin is stable in the sense that the orbits all converge to $x_0$.  However, use of sink orbits as models of stability is called into question with the example of Newhouse: inside his non-hyperbolic family, there are always plane diffeomorphisms with infinitely many sink orbits.

The surface diffeomorphisms studied by Poincar\'e were derived from ODE's and thus not easily computable by direct means.  However, H\'enon showed that chaotic dynamical behaviors arise already from very simple maps, such as $h_{a,b}(x,y)= (1-a x^2 +  y, bx)$.  This family $h_{a,b}$ contains horseshoe maps.  Further, it contains Newhouse examples of:  $(i)$  persistent non-hyperbolicity, and $(ii)$ mappings which have an infinite number of sink basins.   H\'enon performed computer experiments to show that when $a=1.4$ and $b=0.3$, the map $h_{a,b}$ has a ``strange attractor''.  Later, Benedicks-Carleson and others showed that there is a positive measure set of parameters $(a,b)$ for which $h_{a,b}$ actually has an attractor with the properties indicated by the computer in the case of H\'enon's parameters.

Also around the beginning of the 20th century Fatou and Julia began their studies of the iteration of rational functions on the Riemann sphere.  In this case, there is a natural division of dynamical behavior into the stable and chaotic parts, which are carried by what are now called the Fatou and Julia sets.  This elegant theory was developed through extensive use of complex function theory, especially the Montel theory of normal families.  Then, in the early 1980s, this these questions were given a new approach, starting with Douady-Hubbard, Sullivan, Lyubich and others.  This involved the use of the computer, as well as powerful new analytical tools.  In the study of the quadratic family, dynamical bifurcations led to beautiful structures  in parameter space.  

There had also been questions of iteration of (real) polynomial maps of ${\Bbb R}$, strongly influenced by the work of Milnor-Thurston on the combinatorial structure of the dynamics.  Around the same time Feigenbaum and Collet-Tresser introduced the method of renormalization.   If a polynomial has real coefficients, its complexification also yields a map of ${\Bbb C}$.   While these works began strictly within the real domain, the extension to ${\Bbb C}$ permitted the use of complex techniques and has proved very effective.


%

Similarly, the diffeomorphism $h_{a,b}$ extends naturally to a diffeomorphism $H_{a,b}$ of ${\Bbb C}^2$, and we ask about the dynamics of the map $H_{a,b}$.   The virtues of passing to the complex domain in dimensions one and two are similar.  If we can understand $H_{a,b}$, then the real dynamics is obtained by taking the slice with ${\Bbb R}^2$.  If we are looking for periodic points, then we are fixing $n$ and studying the solutions of the polynomial equation pair $h^{n}_{a,b}(x,y)=(x,y)$.  The number of solutions in ${\Bbb R}^2$ may be hard to determine, but the number of solutions (with multiplicity) in ${\Bbb C}^2$ is given by the algebraic degree.  Further, as we vary the coefficients $(a,b)$, the set of solutions will move continuously inside ${\Bbb C}^2$, whereas they may ``appear'' or ``disappear'' from ${\Bbb R}^2$.  In analogy with this, we will expect the dynamics of the complex map $H_{a,b}$ to show more ``completeness''  and ``continuity'' than we find from the real map.  Further, we will have a larger set of mathematical tools that we can bring to bear to the situation.  Besides, the understanding we find in ${\Bbb C}^2$ will be an entry point to real dimension~4.  

In fact, the complex family $H_{a,b}$, where $a$ and $b$ are allowed to be complex, promises to be richer and more interesting than the original real family $h_{a,b}$ discussed above.  Let us consider the family of dynamical systems  given by (generalized) H\'enon maps, which are the polynomial automorphisms of ${\Bbb C}^2$ of the form $f(x,y)=(y,p(y)-\delta x)$, where $p(y)$ is a polynomial of degree $d\ge2$, and $\delta\in{\Bbb C}$, $\delta\ne0$.   (The family $H_{a,b}$ corresponds, modulo linear conjugacy, to  $f(x,y)$ with $p(y)$ quadratic.)  These generate all the polynomial automorphisms of ${\Bbb C}^2$ with positive entropy in the following sense:  Any such  automorphism is (up to conjugacy) a composition of generalized H\'enon maps (see [FM]).   A basic dynamical notion is recurrence. A point is recurrent if iterates of the point come arbitrarily close to it. As a first approximation to this notion we look at the set of bounded orbits. Since $f$ is invertible there are three possible notions of bounded orbits.
 $$K^\pm:=\{ (x,y)\in{\Bbb C}^2: f^{\pm n}(x,y) {\rm\ is\ bounded\ for \ }n\ge0\},\ \ K:=K^+\cap K^-,  \ \ J^\pm:=\partial K^\pm,$$ 
The sets $J^\pm$ represent the boundaries between sets with different behaviors, so we expect to find complicated dynamics at these points.
The sets $J^\pm$  are called the forward/backward Julia sets, and  $J:=J^+\cap J^-$ carries the most chaotic part of the dynamics.  A traditional starting point for the study of this family involves the decomposition into the Julia sets $J^\pm$, and the forward/backward Fatou sets ${\cal F}^\pm:={\Bbb C}^2-J^\pm$, where the behavior is more regular.   A second basic notion is expansion. This means that nearby points are moved apart by some iterate of the map. We can consider forward and backward expansion.  By properties of holomorphic mappings, there is no forward expansion on the interior of $K^+$ and no backward expansion on the interior of $K^-$.

One way to study contracting/expanding properties is to consider stable and unstable sets of a point $p\in{\Bbb C}^2$: 
$$W^{s/u}(p):=\{q\in{\Bbb C}^2:\lim_{n\to+\infty/-\infty}{\rm dist}(f^nq,f^np)=0\}$$
We will discuss some foliations/laminations generated by stable sets:   there are foliations inside the Fatou components, and there are ``laminations'' (in a weak sense) in $J^\pm$, and.  In fact, these ``laminations'' may be very badly behaved and are called ``turbulations''  in [LM].  One can consider the Fatou and Julia sets as independent objects, but here we formulate a number of questions about  the interactions between these sets which define different dynamical regimes.
 
We start with the forward Fatou set ${\cal F}^+$ and its special component $U^+:={\Bbb C}^2-K^+$, consisting of all points which escape to infinity in forward time.   The function 
 $$G^+(x,y):=\lim_{n\to\infty}{1\over d^n}\log^+||f^n(x,y)||$$ 
 gives the (super-exponential) rate of escape to infinity and satisfies $G^+\circ f=d\cdot G^+$.  Further,  $G^+$ is pluri-harmonic on  $U^+=\{G^+>0\}$, so the complex-valued 1-form $\partial G^+$ is holomorphic  and defines a holomorphic foliation ${\cal G}^+$ there.  This foliation is defined in [HO] where it is also shown that the (global) leaves are Riemann surfaces conformally equivalent to ${\Bbb C}$.

 A H\'enon map preserves volume if $|\delta|=1$.  If $|\delta|<1$, then $f$ contracts volume and is said to be {\it dissipative}.  The cases $|\delta|=1$ and $|\delta|\ne1$ are rather different from the point of view of the Fatou components:  When $|\delta|>1$, $K^+$ has no interior, so the whole Fatou set is connected and coincides with $U^+= {\Bbb C}^2-J^+$.  When $|\delta|\le1$,  it is possible for $K^+$ to have interior, but the behavior of $f$ there is different in the cases $|\delta|=1$ and $|\delta|<1$.  When $|\delta|=1$, the interiors of $K^+$ and $K^-$ are the same and coincide with ${\cal F}:={\cal F}^+\cap {\cal F}^-={\rm int}(K)$ (see [FM]).   ${\cal F}$ has finite volume, so any connected component $U$ must be periodic.   A classic result of Newhouse gives a dissipative automorphism $f$ of ${\Bbb R}^2$ which has infinitely many sink orbits.     In the Newhouse example, $f$ can be a H\'enon map with $p(y)$ real and quadratic, and $\delta\ne0$ real and close to zero.  The map $f$ is also a biholomorphic map of ${\Bbb C}^2$, and the complex sink basins of the real sinks give infinitely many components of ${\cal F}^+$.  An important feature of this map is that  it is dissipative.
 
 \proclaim Question 1.  If $f$ preserves volume, can ${\cal F}$ have infinitely many connected components?
 
A H\'enon map with real coefficients gives both a diffeomorphism of ${\Bbb R}^2$ and a diffeomorphism of ${\Bbb C}^2$.  An understanding of the interplay between these two maps, real and complex, will be fundamental.  An area-preserving diffeomorphism of ${\Bbb R}^2$ can have infinitely many ``elliptic islands'' (see, for instance [PD] or [GK]).  The existence of these invariant domains  in ${\Bbb R}^2$ does not answer Question 1 because these elliptic islands are associated with ``twist maps'' and are not contained in Fatou components in ${\Bbb C}^2$, which would  be purely ``rotational'' (see [BS2]).  

If $f$ is a real, area-preserving H\'enon map, and if  $\gamma\subset{\Bbb R}^2$ is a real analytic invariant curve, for instance, a KAM curve, then $\gamma$ complexifies to a maximal Riemann surface $\tilde \gamma\subset{\Bbb C}^2$ which is invariant under $f$ and which is conformally equivalent to an annulus $\{R^{-1}<|\zeta|<R\}$.  We may ask how the real, invariant curves interact with the complex map.  {\it Is it possible for $\tilde \gamma$ to intersect ${\cal F}$?}  If not, then $\tilde \gamma\subset J$.  We let $\tilde \Gamma:=\bigcup \tilde \gamma$ be the union of all complexifications of analytic, invariant curves~$\gamma$.
{\it  If $f$ preserves volume, is $\tilde\Gamma$ dense in $J$?}  We will return to $\tilde \Gamma$ in Question 8.

A more general question about the connection between the (complex) Fatou set on the real map is:

\proclaim Question 2.  Can an area-preserving, real H\'enon map have a Fatou component (other than $U^+$) which intersects ${\Bbb R}^2$?

If $U$ is a periodic component of ${\cal F}$, then we may replace $f$ by $f^N$ and assume that $f(U)=U$.  Continuing with the volume-preserving case, the closure of the iterates $\{f^n|_{U}\}$ generates a compact Abelian group, and its connected component of the identity is either ${\Bbb T}^1$ or ${\Bbb T}^2$ (see [BS2]).  In the case of ${\Bbb T}^2$, $f:U\to U$ is biholomorphically conjugate to a unitary map acting on a Reinhardt domain.  Such a domain can be homeomorphic to either $\Delta^2$, $\Delta\times A$ or $A\times A$, where $\Delta\subset {\Bbb C}$ is the unit disk, and $A=\{1<|\zeta|<R\}$ is a nondegenerate annulus.  The case $\Delta^2$ can occur, and  the case $A\times A$ cannot occur (see [BS2]).  In the case $\Delta^2$, it would be interesting to know more about the boundary of $U$.
  
 \proclaim Question 3.  In the volume-preserving case, can a Fatou component be homeomorphic to $\Delta\times A$?

For the rest of our discussion, we focus on ${\cal F}^+$ in the dissipative case $|\delta|<1$.  Recall that the  components of ${\cal F}^+$ other than $U^+$  are contained in the interior of $K^+$.    (In the dissipative case, we do not consider ${\cal F}^-$.  For a dissipative map, $f^{-1}$ expands volume, so  $K^-(f) = K^+(f^{-1})$, and  $K^-$ has no interior (see [FM]).  The backward Fatou set ${\cal F}^-$ is then connected and  coincides with the  set $U^-$ of points which escape in backward time.)  

A  Fatou component $U$ is said to be {\it wandering} if $f^n(U)\cap U=\emptyset $ for all nonzero $n\in{\Bbb Z}$.  Nonwandering Fatou components are necessarily periodic, i.e.\  $f^N(U)=U$ for some nonzero $N$.  A fundamental question is: 
  
\proclaim Question 4.   {Can there be a wandering Fatou component? }   
  
We recall that in dimension one, a celebrated result of Sullivan shows that a rational map cannot have a wandering domain.  On the other hand, in dimension two, wandering domains have been found for (non-H\'enon) maps in [FS2] and [ABDPR]

Now let $U$ be a periodic component of the interior of $K^+$.  We say that $U$ is {\it recurrent} if there exists $q\in U$ such that $f^n(q)$ does not converge to $\partial U$.  The recurrent components have been classified into two cases: 1.\  basin of an attracting cycle,  2.\  rotational basin  (see [BS2]).  These will be described below.  The remaining (non-recurrent) case is more difficult.  Lyubich and Peters [LP]   show that  if the dissipation satisfies  $ |\delta| < {1\over d^2}$, then a non-recurrent component must be a semi-parabolic basin, which is case 3 below.   It would be important to know whether  the condition on dissipation can be relaxed: 

\proclaim Question 5.   {In the general dissipative case, is every non-recurrent, periodic Fatou component a semi-parabolic basin?}

We briefly describe the three known cases for Fatou components (other than $U^\pm$) in the dissipative case.  (Let us remark, however, that the three cases can coexist, since  different Fatou components can belong to different cases.)

Case 1:  Suppose a H\'enon map $f$ has an attracting fixed point $O$, and suppose that the multipliers at $O$ are $0<|\alpha|<|\beta|<1$.  We consider the basin of attraction (stable set) ${\cal B}:=W^s(O)$.  If there is no resonance between $\alpha$ and $\beta$, then the restricted map $(f,{\cal B})$ is conjugate to the linear map $(L,{\Bbb C}^2)$, where $L(z,w)= (\alpha z,\beta w)$.   Whether there is a resonance or not, there is a holomorphic map $\Phi:{\cal B}\to {\Bbb C}$ such that $\Phi\circ f = \beta\cdot \Phi$.   (If $f$ is linearizable, then $\Phi$ is one of the coordinates of the linearizing conjugacy.)   The strong stable manifold  of a point $p\in{\cal B}$ is defined as the points which converge as rapidly as possible:
 $$W^{ss}(p)=\{q\in{\cal B}: \lim_{n\to\infty}{1\over n}\log ({\rm dist}(f^n(q),f^n(p)))=\log|\alpha|\}$$
Let ${\cal W}^{ss}({\cal B})$ denote the strong stable foliation, which is generated by the sets $W^{ss}(p)$, $p\in{\cal B}$.  It follows that ${\cal W}^{ss}$ consists of the level sets of $\Phi$, and thus is also a fibration by Riemann surfaces which are closed in ${\cal B}$.

Case 2.  In this case we suppose that the Fatou component $U$ contains an invariant Riemann surface ${\cal M}\subset{\Bbb C}^2$, and the restriction $f|_{\cal M}$ acts as a rotation.  There are two cases:  ${\cal M}$ can be conformally equivalent to a disk or an annulus.  These are reminiscent of the one-dimensional cases of  a Siegel disk or a Herman ring.  An important open question for H\'enon maps is:  

\proclaim Question 6.  {Can the case of the annulus actually occur?}  

Let ${\cal R}\subset{\Bbb C}$ denote the disk $\Delta$ or an annulus $\{1<|\zeta|<R\}$, depending on the case.  Then there is a number $\omega$, $|\omega|=1$, $\omega$ not a root of unity, and a biholomorphic map $\Phi:{\cal R}\times{\Bbb C}\to U$ which conjugates the linear map $L(z,w) = (\omega z, \delta w/\omega)$ to $f|_U$ (see [BS2]).  Recall that for a set $X$, we define its stable set 
$$W^s(X):= \{q\in{\Bbb C}^2: \lim_{n\to\infty}{\rm dist}(f^nq,f^n(X))=0\}.$$ 
Thus $U=W^s({\cal M})$ is the basin of attraction of ${\cal M}$.  Further, for $\zeta\in{\cal R}$, $\Phi(\{\zeta\}\times{\Bbb C})= W^s(\Phi(\zeta,0))$, so we have a foliation (and fibration) ${\cal W}^s$ of $U$ by stable manifolds.  

Case 3.   A semi-attracting/semi-parabolic point is a fixed (or periodic) point $O$ with multipliers $1$ and $\delta$ with $|\delta|<1$.  There is a semi-parabolic basin ${\cal B}$ for $O$ and an Abel-Fatou function $\Phi_{AF}:{\cal B}\to {\Bbb C}$, which satisfies $\Phi_{AF}\circ f=\Phi_{AF}+1$ (see [U]).  There is a strong stable foliation ${\cal W}^{ss}({\cal B})$ as before, and the strong stable manifolds have been shown to be the sets where $\Phi_{AF}$ is constant and are conformally equivalent to ${\Bbb C}$  (see [BSU]).  

For the Fatou components $U$ in any of the cases 1, 2 and 3 discussed above,  we may let $M$ denote any of the leaves of the invariant foliations.  While $M$ is a closed subset of ${\cal B}$, it is not closed in ${\Bbb C}^2$, and  by [BS2] its ${\Bbb C}^2$-closure,  $\overline M$, contains $J^+$.  As a consequence, $\partial U=J^+$.  

The sets $J^\pm$ have no interior and thus cannot carry foliations.  However, they contain many Riemann surfaces:  For every saddle point $p$, the unstable manifold $W^u(p)\subset J^+$, and $W^u(p)$ is dense in $J^+$.  There are many such unstable manifolds, since by [BLS], the number of saddle points of period $n$ is asymptotically $d^n$ as $n\to\infty$.  

Saddle points show both recurrence and expansion.  The set $S$ of saddle points is very important for understanding $f$.  We will denote its closure by $J^*$.  The set $J^*$ has several dynamical characterizations, including being the support of the unique measure of maximal entropy (see [BLS]).  Thus $J^*$ is a strong candidate for being called the Julia set, and it is fundamental to ask: 

\proclaim Question 7.  Is $J=J^*$?

$J^*$ is always contained in $J$.  In case $J\ne J^*$, it would be important to know the nature of the set $J-J^*$.  For instance:

\proclaim Question 8.  If \  $\tilde\Gamma$ is the set defined after Question 2, is $\tilde\Gamma\subset J^*$? 

We say that a set ${\cal L}$ carries a Riemann surface lamination if  it looks locally like a product $\Delta\times A$, where  $\Delta\subset{\Bbb C}_x$ is the unit disk, and $A\subset{\Bbb C}_y$ is compact.  More precisely,  there is a local homeomorphism in $(x,y)$ which is holomorphic in $x$.  Every foliation is also a lamination, and conversely, laminations may be thought of as foliations of compact sets.   

We say that a H\'enon map $f$ is hyperbolic if $J$ is a hyperbolic set, which means that there are stable/unstable manifolds ${\cal W}^{s/u}=\{W^{s/u}(x):x\in J\}$ which give  Riemann surface laminations of $J^\pm$ in a neighborhood of $J$.  (And conversely, it was shown in [BS8] that the existence of transverse laminations of $J^\pm$ gives hyperbolicity.)  Each $W^{s/u}(x)$ is conformally equivalent to ${\Bbb C}$.   If a H\'enon map is hyperbolic, then the interior of  $K^+$ consists of the basins of finitely many sink orbits.  Further, $J$ coincides with the non wandering set, $J=J^*$, and $f|_{J}$ is structurally stable (see [BS1] for a treatment of the hyperbolic case).

While hyperbolic is the class of maps that is best understood, hyperbolic maps are rather special.  Hyperbolicity requires the simultaneous conditions of uniform expansion and contraction.  Below, we will discuss the concepts of quasi-expansion and quasi-contraction which allow us to deal with more general maps that have weaker sorts of expansion/contraction; and it allows us to consider expansion and contraction separately.

We remark that in general (even in the absence of hyperbolicity) there is a measure-theoretic sense in which $J^+$ is filled a.e.\   by Riemann surfaces.  The reason for this is  that $J^+$ carries a unique invariant ``laminar'' current, and this is constructed from stable manifolds (see [BLS], [D1--3]).  In fact, there is even a stronger (non dynamical)  sense of uniqueness:  the only positive closed currents supported in $K^+$ are multiplies of $\mu^+$ (see [FS1], [DS]).   Despite the central importance of currents, we do not discuss them here.

In the hyperbolic case, we have two  families of Riemann surfaces:  the foliation ${\cal G}^+$ of $U^+$, and the lamination ${\cal W}^s(J)$ of $J^+$.  We can ask how these families interact.  For instance, do they fit together continuously?  One formulation of ``continuity to the boundary'' for ${\cal G}^+$ is the condition that ${\cal G}^+\cup {\cal W}^s(J)$ is itself a lamination.  This was shown in [BS7] to be the case when $J$ is connected.  

However, if $J$ is hyperbolic but not connected, then  ${\cal G}^+\cup {\cal W}^s(J)$  is not a lamination at points of $J$.  To see this, let  $q\in J$ be a saddle point, and let $W^u(q)$ denote the unstable manifold of $q$.  Let ${\cal T}_q$ denote the set of tangencies between ${\cal G}^+$ and $W^u(q)$.  By [BS6], ${\cal T}_q\ne\emptyset$ whenever $J$ is not connected.  On the other hand, $W^u(q)$ is transverse to ${\cal W}^s(J)$, which means that cannot have the local product structure necessary to be a lamination at $q$.  And when $J$ is hyperbolic, the saddle points are dense in $J$.  By the Lambda Lemma, it follows that  ${\cal G}^+\cup {\cal W}^s(J)$ cannot be laminar anywhere at $J^+$.

We have a similar situation in Case 1.  In [DL] it is shown that if $ |\delta| < {1\over d^4}$, there are tangencies between  $W^u(p)$ and ${\cal W}^{ss}({\cal B})$.  It follows that in the hyperbolic case, ${\cal W}^{ss}({\cal B})\cup {\cal W}^s(J)$ is not a lamination.  [DL] also shows that if the dissipation satisfies $ |\delta| < {1\over d^2}$, there are tangencies in case 3.  However, in cases 2 and 3, $f$ cannot be hyperbolic.

%
%
\proclaim Question 9.  When is it possible for $J^+$ to carry a Riemann surface lamination ${\cal L}^+$?

By Radu and Tanase [RT], the answer is ``yes'' in the semi-parabolic case 3 if $f$ is quadratic and $|a|\ll 1$. {\it  Is the answer always ``no'' in the rotational  case 2?}

%

Let ${\cal M}$ be a family of complex manifolds $M\subset {\Bbb C}^2$.  We  introduce a condition on ${\cal M}$ that is more inclusive than the condition of laminarity.     Given a point $s$, we let $B(s,\epsilon)$ denote the $\epsilon$-ball centered at $s$.  We say that the family ${\cal M}$ is {\it locally proper at the set $S$, in the strong sense,} if there is an $\epsilon>0$ such that for each $s\in S$ and each $M\in{\cal M}$, every connected component $M'$ of $M\cap B(s,\epsilon)$ is closed in $B(s,\epsilon)$.  Thus each component $M'$ is a subvariety of $B(s,\epsilon)$.  We say that  the family ${\cal M}$ is {\it locally proper at the set $S$, in the weak sense, } if there is an $\epsilon>0$ such that for each $s\in \bar S$ and each $M\in{\cal M}$ containing~$s$, $M_s$ is closed in $B(s,\epsilon)$, where $M_s$ denotes the connected component of $M\cap B(s,\epsilon)$ which contains $s$.

We say that ${\cal M}$ has {\it locally  bounded area  at the set $S$, in the strong sense } if  the supremum over all components $M'$ of $M\cap B(s,\epsilon)$ satisfies  $\sup_{M'}Area(M')<\infty $.  The weak version is: ${\cal M}$ has {\it locally bounded area  at $S$, in the weak sense } if $\sup_{s\in S}Area(M_s)<\infty$, where  $M_s$ denotes the component of $M\cap B(s,\epsilon)$ containing~$s$.

{\it Quasi-expansion} a loosening of the notion of uniform expansion.  In one dimension, it is closely related to the very useful Misiurewicz property and semi-hyperbolicity (see [CJY]).   Several equivalent definitions of quasi-expansion are given in [BS8].  Some of them involve uniform expansion in terms of (discontinuous) metrics which are not necessarily equivalent to the Euclidean metric.  One of them involves the locally proper bounded area condition  in terms of the set $S=\{\rm saddle\ points\}$: $f$ is quasi-expanding if and only if $\{W^u(p):p\in S\}$ is locally proper with locally bounded area, in the weak sense.  A consequence (see [BS8]) is that if $f$ is quasi-expanding, then there are unstable manifolds through each point of $\bar S$, although the unstable manifolds do not necessarily have the local product structure of a lamination.  The obstruction to local product structure is local folding on arbitrarily small scales.  Although local product structure may fail, the area bound limits the degree of local folding.  If ${\cal M}=\{W^u(p):p\in S\}$ is locally proper with locally bounded area, in the strong sense, then these stable manifolds fill out $J^-$ in a neighborhood of $J^*$. 

We say that $f$ is {\it quasi-contracting} if $f^{-1}$ is quasi-expanding.  In Case 3, $f$ can not be quasi-expanding, but it is not known when it can be quasi-contracting.  In Case 2, it is not known whether $f$ can be quasi-expanding or quasi-contracting.  
We ask about boundary behaviors of the various  foliations that arise: we do not expect ``continuity'' (i.e. the union of the two families is again laminar) but we can hope for some bound on the amount of folding:
\proclaim Question 10.  If $f$ is hyperbolic (or merely quasi-contracting), does the foliation ${\cal G}^+$ have  locally proper bounded area at points of $J^+$? 


Similarly, we ask the same question for laminations of the Fatou components of ${\rm int}(K^+)$:

\proclaim Question 11.  {In cases 1, 2 or 3, when can ${\cal W}^{ss}({\cal B})$ have locally  proper, bounded area at points of $J^+$?} 

In the dissipative case with $J$ connected, it was shown in [BS6] that $J^--K$ carries a Riemann surface lamination, which we will call ${\cal J}^-$.  In the hyperbolic, dissipative case, ${\cal G}^-\cup{\cal J}^-$ is not a Riemann surface lamination.

\proclaim Question 12.  In the dissipative, $J$ connected case,  does ${\cal J}^-$ have the proper, locally bounded area condition?  What if, in addition, $J$ is hyperbolic?

\noindent{\it Acknowledgment }  I wish to thank S. Cantat, R. Dujardin, M. Lyubich, and J. Smillie for helpful comments on this essay.

\bigskip
 \centerline{\bf References}
 
\medskip
\item{[ABDPR]} M. Astorg, X. Buff, R. Dujardin, H. Peters, J. Raissy,   A two-dimensional polynomial mapping with a wandering Fatou component.  {\tt  arXiv:1411.1188}

\item{[BLS]}  Eric Bedford,  Mikhail Lyubich and John Smillie,    Polynomial diffeomorphisms of ${\bf C}^2$. IV. The measure of maximal entropy and laminar currents. Invent. Math. 112 (1993), no. 1, 77--125.
 
\item{[BS1]}   Eric Bedford and John Smillie, Polynomial diffeomorphisms of ${\bf C}^2$:   Currents, equilibrium measure and hyperbolicity. Invent. Math. 103 (1991), no. 1, 69--99. 
 
\item{[BS2]} Eric Bedford and John Smillie, Polynomial diffeomorphisms of ${\bf C}^2$. II. Stable manifolds and recurrence. J. Amer. Math. Soc. 4 (1991), no. 4, 657--679.

\item{[BS6]}  Eric Bedford and John Smillie, Polynomial diffeomorphisms of ${\bf C}\sp 2$. VI. Connectivity of $J$. Ann. of Math. (2) 148 (1998), no. 2, 695--735.

\item{[BS7]}  Eric Bedford and John Smillie, Polynomial diffeomorphisms of ${\bf C}^2$. VII. Hyperbolicity and external rays. Ann.\ Sci.\ \'Ecole Norm.\ Sup.\ (4) 32 (1999), no. 4, 455--497.

\item{[BS8]}  Eric Bedford and John Smillie, Polynomial diffeomorphisms of ${\bf C}^2$. VIII. Quasi-expansion. Amer. J. Math. 124 (2002), no. 2, 221--271.

 \item{[BSU]}  Eric Bedford, John Smillie and Tetsuo Ueda,  Semi-parabolic bifurcations in complex dimension two.  {\tt  arXiv:1208.2577 }
 
\item{[CJY]}  L. Carleson, P. Jones and J-C Yoccoz,   Julia and John. Bol. Soc. Brasil. Mat. (N.S.) 25 (1994), no. 1, 1--30. 
 
\item{[DS]}  T-C Dinh and N. Sibony,   Rigidity of Julia sets for H\'enon type maps, {\tt arXiv:1301.3917}
 
\item{[PD]}  P. Duarte,  Elliptic isles in families of area-preserving maps. Ergodic Theory Dynam. Systems 28 (2008), no. 6, 1781--1813. 
 
\item{[D1]}  R. Dujardin, Structure properties of laminar currents on ${\bf  P}\sp 2$. J. Geom.\ Anal.\ 15 (2005), no.\ 1, 25--47.

\item{[D2]} R. Dujardin, Sur l'intersection des courants laminaires.   Publ. Mat. 48 (2004), no. 1, 107--125.

\item{[D3]} R. Dujardin, Laminar currents in ${\bf P}\sp 2$. Math. Ann. 325 (2003), no. 4, 745--765.
 
\item{[DL]} R. Dujardin and M. Lyubich, Stability and bifurcation for dissipative polynomial automorphisms of ${\bf C}^2$, Inventiones math.  2014. {\tt arXiv:1305.2898}

\item{[FS1]}  J. E. Forn\ae ss and N. Sibony,  Complex dynamics in higher dimensions, Complex
Potential Theory, ed.\ by P. M. Gauthier. Kluwer Academic Publishers, 1994.

\item{[FS2]}  J. E. Forn\ae ss and N. Sibony,  Fatou and Julia sets for entire mappings in ${\bf C}^k$.
Math. Ann. 311 (1998), no. 1, 27--40. 

\item{[FM]}  Shmuel Friedland and John Milnor,
Dynamical properties of plane polynomial automorphisms.
Ergodic Theory Dynam. Systems 9 (1989), no. 1, 67--99. 

\item{[GK]}  A. Gorodetski and V. Kaloshin,  Conservative homoclinic bifurcations and some applications,  Proceedings of the Steklov Institute of Mathematics, December 2009, Volume 267, Issue 1, pp 76--90 

\item{[HO]}  Hubbard, J. H.; Oberste-Vorth, R. W. H\'enon mappings in the complex domain. I. The global topology of dynamical space. IHES Sci. Publ. Math. No. 79 (1994), 5--46.

\item{[LM]}  Lyubich, Mikhail and Minsky, Yair, Laminations in holomorphic dynamics.
J. Differential Geom. 47 (1997), no. 1, 17--94. 

\item{[LP]}  Lyubich, Mikhail and Peters, Han,  Classification of invariant Fatou components for dissipative H\'enon maps, Geom.\ Funct.\ Anal.\ Vol. 24 (2014) 887--915.

 \item{[RT]}  R. Radu and R. Tanase,  A structure theorem for semi-parabolic H\'enon maps.  
 
 {\tt  arXiv:1411.3824}
 
%
\item{[U]}   Tetsuo Ueda, Local structure of analytic transformations of two complex variables. I.  J. Math. Kyoto Univ. 26 (1986), no. 2, 233--261.

\bigskip\bigskip 

 \rightline{Stony Brook University}
 
 \rightline {Stony Brook, NY 11794}
 
 \rightline{\tt ebedford@math.sunysb.edu}
 
\bye

However, if $J$ is not connected, then  ${\cal G}^+\cup {\cal W}^s$  is not a lamination at points of $J$.  To see this, let  $q\in J$ be a saddle point, and let $W^u(q)$ denote the unstable manifold of $q$.  Let ${\cal T}_q$ denote the set of tangencies between ${\cal G}^+$ and $W^u(q)$.  By [BS6], ${\cal T}_q\ne\emptyset$ whenever $J$ is not connected.  On the other hand, $W^u(q)$ is transverse to ${\cal W}^s$, which means that  ${\cal G}^+\cup {\cal W}^s$ cannot be a lamination at $q$.  When $J$ is hyperbolic, the saddle points are dense in $J$.  The question remains:  {\it  Is  ${\cal G}^+\cup {\cal W}^s$ a lamination at points of $J^+-J$?}

Let us consider a condition which is more inclusive.   Let $\Gamma$ denote a family of Riemann surfaces in ${\Bbb C}^2$.  We say that $\Gamma$ satisfies the {\it locally proper, bounded area} condition at $J^+$ if for each $p\in J^+$, there is a neighborhood $U$ of $p$ such that for each surface $S\in \Gamma$, each connected component $S'$ of $U\cap S$ is closed in $U$, and $Area(S')$ is uniformly bounded.  Let us recall that a Riemann surface lamination may be written locally in a local neighborhood $U$ which is equivalent to the bidisk $\Delta\times\Delta$ as a union of disjoint holomorphic graphs $\gamma_a:= \{y=\varphi_a(x): x\in \Delta\}$, with $a\in A$.  The condition that $\gamma_a$ is a graph means that it is a subvariety of $\Delta_x\times\Delta_y$ which has a bijective projection to $\Delta_x$.   Locally proper, bounded area is equivalent to requiring that there is a local neighborhood $U$ equivalent to a bidisk, and ${\cal L}\cap U$ is given by $\bigcup_{a\in A}\gamma_a$, where $\gamma_a$ is a subvariety of $U$, and the projection of $\gamma_a$ to $\Delta_x$ is proper with degree bounded independently of $a$.   

Let us remark that the condition of ``locally proper bounded area'' arose in [BS8] where we define ``quasi-expansion''.  Let ${\cal S}$ denote the set of saddle points for $f$, and for each $p\in{\cal S}$ we let $W^u(p)$ denote the unstable manifold.  Then $\bigcup_{p\in {\cal S}}W^u(p)$ has the locally proper bounded area condition if and only if $f$ is quasi-expanding.  The map $f$ is said to be {\it quasi-hyperbolic} if both $f$ and $f^{-1}$ are quasi-expanding.

Returning to the question above, we can ask: {\it If $f$ is hyperbolic, does ${\cal G}^+\cup {\cal W}^s$ satisfy the locally proper, bounded area condition?}  More generally, we can ask whether this holds for quasi-hyperbolic maps.

Suppose a H\'enon map $f$ has an attracting fixed point $O$, and suppose that the multipliers at $O$ are $0<|\alpha|<|\beta|<1$.   If there is no resonance between $\alpha$ and $\beta$, then the restricted map $(f,{\cal B})$ is conjugate to the linear map $(L,{\Bbb C}^2)$, where $L(z,w)= (\alpha z,\beta w)$.   Whether there is a resonance or not, there is a holomorphic map $\Phi:{\cal B}\to {\Bbb C}$ such that $\Phi\circ f = \beta\cdot \Phi$.   The strong stable manifold  of a point $p\in{\cal B}$ is defined as:
 $$W^{ss}(p)=\{q\in{\cal B}: \lim_{n\to\infty}{1\over n}\log ({\rm dist}(f^n(q),f^n(p)))=\log|\alpha|\}$$
Let ${\cal W}^{ss}$ denote the strong stable foliation, which is generated by the sets $W^{ss}(p)$, $p\in{\cal B}$.  It follows that ${\cal W}^{ss}$ consists of the level sets of $\Phi$, and thus is also a fibration.
If $J$ is hyperbolic, we can ask:  {\it  Does  ${\cal W}^{ss}({\cal B})\cup{\cal W}^s(J)$ have the locally proper, bounded area property?}

We note that there are cases where  ${\cal W}^{ss}({\cal B}) \cup{\cal W}^s(J)$ cannot be a lamination.  We say that $f$ is {\it moderately dissipative} if the complex Jacobian determinant of $f$ is $<{1\over d}$ ($d={\rm degree}(f)$.   It has been shown in [DL] that if $f$ is moderately dissipative, then ${\cal T}_q\ne\emptyset$.  
 
 The situation is analogous for maps with a semi-attracting/semi-parabolic point, i.e.\  a fixed point $O$ with multipliers $1$ and $a$ with $|a|<1$.  There is a semi-parabolic basin ${\cal B}$ for $O$.  There is an Abel-Fatou function $\Phi_{AF}:{\cal B}\to {\Bbb C}$, which satisfies $\Phi_{AF}\circ f=\Phi_{AF}+1$.  There is a strong stable foliation ${\cal W}^{ss}({\cal B})$ as before, and the strong stable manifolds have been shown to be the sets where $\Phi_{AF}$ is constant.  
%
 We can ask:  {\it Is it possible for $J^+$ to carry a Riemann surface lamination ${\cal L}^+$?}  (By Radu and Tanase, the answer is ``yes'' if $f$ is quadratic and $|a|\ll 1$.)  In this case, we may ask:  {\it Is ${\cal W}^{ss}({\cal B})\cup{\cal L}^+$ or ${\cal G}^+\cup {\cal L}^+$ a Riemann surface lamination in a neighborhood of $J^+-J$?}  For a saddle $r$, we  define the set ${\cal T}_r$ of tangencies  between  ${\cal W}^{ss}({\cal B})$ and $W^u(r)$.  Dujardin and Lyubich showed that if $f$ is moderately dissipative, then ${\cal T}_r\ne\emptyset$, so the laminarity of  ${\cal W}^{ss}({\cal B})\cup{\cal L}^+$ cannot hold on $J$ itself. 
 
 Finally, we can also ask about the case of Siegel basins.

 \bigskip\bigskip

\item{[BLS]}  Bedford, E.; Lyubich, M.; Smillie, J. Distribution of periodic points of polynomial diffeomorphisms of $\bf C\sp 2$. Invent. Math. 114 (1993), no. 2, 277--288.
 
\item{[BS2]} Eric Bedford and John Smillie, Polynomial diffeomorphisms of ${\bf C}^2$. II. Stable manifolds and recurrence. J. Amer. Math. Soc. 4 (1991), no. 4, 657--679.

\item{[BS6]}  Eric Bedford and John Smillie, Polynomial diffeomorphisms of ${\bf C}\sp 2$. VI. Connectivity of $J$. Ann. of Math. (2) 148 (1998), no. 2, 695--735.

\item{[BS7]}  Eric Bedford and John Smillie, Polynomial diffeomorphisms of ${\bf C}^2$. VII. Hyperbolicity and external rays. Ann.\ Sci.\ \'Ecole Norm.\ Sup.\ (4) 32 (1999), no. 4, 455--497.

\item{[BS8]}  Eric Bedford and John Smillie, Polynomial diffeomorphisms of ${\bf C}^2$. VIII. Quasi-expansion. Amer. J. Math. 124 (2002), no. 2, 221--271.

 \item{[BSU]}  Eric Bedford, John Smillie and Tetsuo Ueda,  Semi-parabolic bifurcations in complex dimension two.  
 
\item{[D1]}  R. Dujardin, Structure properties of laminar currents on ${\bf  P}\sp 2$. J. Geom.\ Anal.\ 15 (2005), no.\ 1, 25--47.

\item{[D2]} R. Dujardin, Sur l'intersection des courants laminaires.   Publ. Mat. 48 (2004), no. 1, 107--125.

\item{[D3]} R. Dujardin, Laminar currents in ${\bf P}\sp 2$. Math. Ann. 325 (2003), no. 4, 745--765.
 
\item{[DL]} R. Dujardin and M. Lyubich, Stability and bifurcation for dissipative polynomial automorphisms of ${\bf C}^2$, Inventiones math.  2014. {\tt arXiv:1305.2898}

\item{[FM]}  Shmuel Friedland and John Milnor,  Dynamical properties of plane polynomial automorphisms.
Ergodic Theory Dynam.\ Systems 9 (1989), no. 1, 67--99. 

\item{[HOV]}  Hubbard, J. H.; Oberste-Vorth, R. W. H\'enon mappings in the complex domain. I. The global topology of dynamical space. IHES Sci. Publ. Math. No. 79 (1994), 5--46.

\item{[LP]}  Lyubich, Mikhail and Peters, Han,  Classification of invariant Fatou components for dissipative H\'enon maps, Geom.\ Funct.\ Anal.\ Vol. 24 (2014) 887--915.

 \item{[RT]}  R. Radu and R. Tanase,  A structure theorem for semi-parabolic H\'enon maps.
 
 \bye